\input amstex
\input xypic
\documentstyle{amsppt}
\input psfig
\def \sm{\sigma}
\magnification=1200
\expandafter\redefine\csname logo\string@\endcsname{}
\NoBlackBoxes
\NoRunningHeads

\topmatter
\title A representation of generalized braid group in classical braid
group\\
\endtitle
\date School of Mathematics, Tata Institute, India. \enddate
\author S. K. Roushon
\endauthor
\address School of Mathematics, Tata Institute of Fundamental 
Research, Homi Bhabha Road, Mumbai 400 005, India. 
\endaddress
\email roushon\@math.tifr.res.in \endemail
\thanks e-mail:- roushon\@math.tifr.res.in and fax:- +91 22 215 2110/2181  
\endthanks
\abstract
We ask if any finite type generalized braid group is a subgroup of some 
classical Artin braid group. We define a natural map from a given finite 
type generalized braid group to a classical Artin braid group and ask if
this map is an injective homomorphism. We also check that this map is 
a homomorphism for the generalized braid groups of type $A_n$, $B_n$ 
and $I_2(k)$.

\endabstract
\endtopmatter
\document
\baselineskip 14pt

\head
0. Introduction
\endhead
The main impetus to this article came from the paper \cite{8} where 
it was proved that the Whitehead group of any subgroup of the classical 
(Artin) braid group vanish. A generalization of the Artin braid group 
is the generalized braid group associated to Coxeter groups. This 
important class of groups are not very well understood. Naturally one 
asks if the similar theorem as in \cite{8} is true for generalized 
braid groups. We consider only finite type generalized braid groups. 
These groups can be presented in terms of generators and relations, from
finite Coxeter groups after removing some of the relations in the Coxeter
group presentation. We recall that any finite type generalized 
braid group is torsion free. In fact these groups are fundamental groups  
of finite aspherical complexes (\cite{7}). 

In fact in \cite{8} it was proved that the Whitehead group of a torsion 
free subgroup of any finite extension of a strongly poly-free group 
vanish. (The class of strongly poly-free groups were introduced in 
\cite{1}). In \cite{1} it was proved that the classical pure 
braid group is strongly poly-free. On the other hand it is not yet known 
if the pure braid group associated to a generalized braid group 
is strongly poly-free. It is not even known if these groups are 
poly-free (recently this question is also asked by M. Bestvina in 
\cite{2}). Though for certain class of Coxeter groups it is 
known that the associated pure braid group is polyfree (see \cite{4}). 
The pure braid group associated to a generalized 
braid group is the kernel of the obvious homomorphism from the 
generalized braid group onto the associated Coxeter group. A 
complete classification of finite Coxeter groups is known 
(see Table 1) (\cite{6}). Here recall that the symmetric 
groups are Coxeter group of type $A_n$ and the associated 
generalized braid group is the classical Artin braid group. 

As any finite group can be embedded in a symmetric group, 
we embed a finite Coxeter group in a symmetric group and 
extend this embedding to a map from the associated generalized 
braid group to the classical braid group. We ask if we can 
choose this extension to be an injective homomorphism. 

In this article we check that we can choose this extension 
to be a homomorphism for the generalized braid group of type 
$A_n$, $B_n$ and $I_2(k)$. Whether this homomorphism is injective 
is not yet settled. 

In Section 4 we give examples of some natural extension in the 
case of braid group of type $D_n$ which is not a homomorphism. 

We prove the following theorem which gives a partial answer to our
question:

\proclaim{Theorem} Let ${\Cal W}_S$ be a Coxeter group of type 
$A_n$, $B_n$ or $I_2(k)$, then there is an embedding $e:{\Cal W}_S\to 
S_m$ for some $m$ and a homomorphism ${\Cal A}_S\to {\Cal A}_m$ making 
the following diagram commutative:\endproclaim

$$\diagram
1\rto  & {\Cal {PA}}_S\rto \dto^{{f_S}|_{{\Cal
{PA}}_S}} & {\Cal A}_S
\rto^{\pi_S} \dto^{f_S} & {\Cal W}_S
\rto \dto^e & 1 \\   
1\rto & {\Cal {PA}}_m\rto & {\Cal A}_m \rto^{\pi_m} &  S_m \rto & 1
\enddiagram
$$

Here ${\Cal A}_S$ is the braid group associated to the Coxeter 
group ${\Cal W}_S$ and ${\Cal A}_m$ is the classical Artin braid 
group on $(m-1)$-strings.

This paper is an expanded and rewritten version of \cite{11}.

\head
1. The question and some consequences
\endhead

Let $S=\{s_1,s_2,\ldots,s_k\}$ be a finite set and $m:S\times S\to
\{1,2,\cdots, \infty\}$ be a 
function with the property $m(s,s)=1$ and $m(s,s')\geq 2$ for $s\neq s'$. 
The Coxeter group associated to the system $(S,m)$ is by definition the 
group ${\Cal W}_S=\{S\ |\ (ss')^{m(s,s')}=1,\ s,s'\in S\}$ and no relation
if $m(s,s')=\infty$. Throughout we always assume that the group ${\Cal
W}_S$ is finite. In this special case $m$ is always finite. A complete 
classification of finite Coxeter groups are known (see \cite{6}). We 
reproduce the enumeration of all finite Coxeter groups in Table 1. The 
symmetric groups $S_n$ on $n$ letters are examples of Coxeter groups. 
These are the Coxeter group of type $A_n$ in the table. 
The generalized braid group associated to the Coxeter
group ${\Cal W}_S$ is ${\Cal A}_S=\{S\ |\ ss'ss'\cdots = s'ss's\cdots,
s,s'\in S\}$,
here the number of times the factors in $ss'ss'\cdots$ appears is 
$m(s,s')$, i.e., if $m(s,s')=3$ then the relation is $ss's=s'ss'$. In 
such a case we call ${\Cal A}_S$ a braid group of type 
${\Cal W}_S$. There is an obvious surjective homomorphism 
${\Cal A}_S\to {\Cal W}_S$. Throughout this 
paper by {\it classical braid group} we will mean the generalized 
braid group associated to the Coxeter group $A_n$. 

Suppose the group ${\Cal W}_S$ is finite. Consider the 
symmetric group $S_m$ on $m$ letters. It is generated by $N=m-1$ elements 
of order $2$. We have a surjective homomorphism ${\Cal A}_m\to S_m$, 
where ${\Cal A}_m$ is the classical Artin braid group on $(m-1)$-strings. 
We choose an embedding $e:{\Cal W}_S\to S_m$ (for example let $m\geq 
|{\Cal W}_S|$). And this embedding can be extended to give a  
map $f_S:{\Cal A}_S\to {\Cal A}_m$ defined on the generators 
of ${\Cal A}_S$ in the same way the map ${\Cal W}_S\to S_m$ is given 
such that the following diagram commutes.

$$\diagram
1\rto  & {\Cal {PA}}_S\rto \dto^{{f_S}|_{{\Cal
{PA}}_S}} & {\Cal A}_S
\rto^{\pi_S} \dto^{f_S} & {\Cal W}_S
\rto \dto^e & 1 \\
1\rto & {\Cal {PA}}_m\rto & {\Cal A}_m \rto^{\pi_m} &  S_m \rto & 1 
\enddiagram
$$

We call this diagram of groups as diagram $D$.

Here ${\Cal PA}_S$ and ${\Cal PA}_m$ stand for the pure braid groups of
the groups ${\Cal A}_S$ and ${\Cal A}_m$ respectively. 

Our question is that:

\proclaim{Question 1} Is the map $f_S$ defined above a homomorphism? If
yes, then is it injective? If it is not injective then what is the kernel 
of $f_S$?\endproclaim

Note that injectivity of $e$ implies that ker$(f_S)$ is same as 
ker$({f_S}|_{{\Cal PA}_S})$. As the pure braid group is the  
fundamental group of the complement of reflection hyperplanes 
in the complex $n$-space corresponding to a faithful representation 
in $GL(n, {\Bbb R})$ of the associated Coxeter group as a relection 
group, some geometric method might be helpful to prove the injectivity 
of the map $f_S$.

We describe the map $f_S$ in little more details: suppose $S_m$ is 
generated by the symbols $r_1,r_2,\ldots,r_{m-1}$ and the relations are 
$r_ir_j=r_jr_i$ whenever $|i-j|\geq 2$, 
$r_ir_{i+1}r_i=r_{i+1}r_ir_{i+1}$ for $i\leq m-2$ and $r_i^2=1$. 
Then ${\Cal A}_m=\{r_1,r_2,\ldots,r_N\ |\ r_ir_j=r_jr_i\ \text{whenever}\ 
|i-j|\geq 2, r_ir_{i+1}r_i=r_{i+1}r_ir_{i+1}\ \text{for}\  i\leq m-2\}$.
The surjective homomorphism ${\Cal A}_m\to S_m$ is the obvious map 
defined on the generators.  

Now define the map $f_S$ by: $f_S(s_i)$ is a coset representative of the
coset containing $e(s_i)$ when $e(s_i)$ is 
considered as an element in ${\Cal A}_m$. Here we note that the 
answer to the question depends on the choice of the coset
representative. In the next section we make a suitable 
choice of the coset representative, in the cases of the generalized 
braid groups of type $A_n$, $B_n$ and $I_2(k)$, to show
that the map $f_S$ is a homomorphism.

If the answers to the first two questions above are yes then a
consequence of the main theorem in \cite{8} will be:

\proclaim{Consequence 1.1} The Whitehead group, reduced projective class 
group and the lower $K$-groups of any subgroup 
of a finite type generalized braid group vanish.\endproclaim

Also another useful consequence will be:

\proclaim{Consequence 1.2} If the classical Artin braid group has a 
discrete faithful linear representation then the same is true for any 
finite type generalized braid group. \endproclaim 

\medskip

\centerline{\psfig{figure=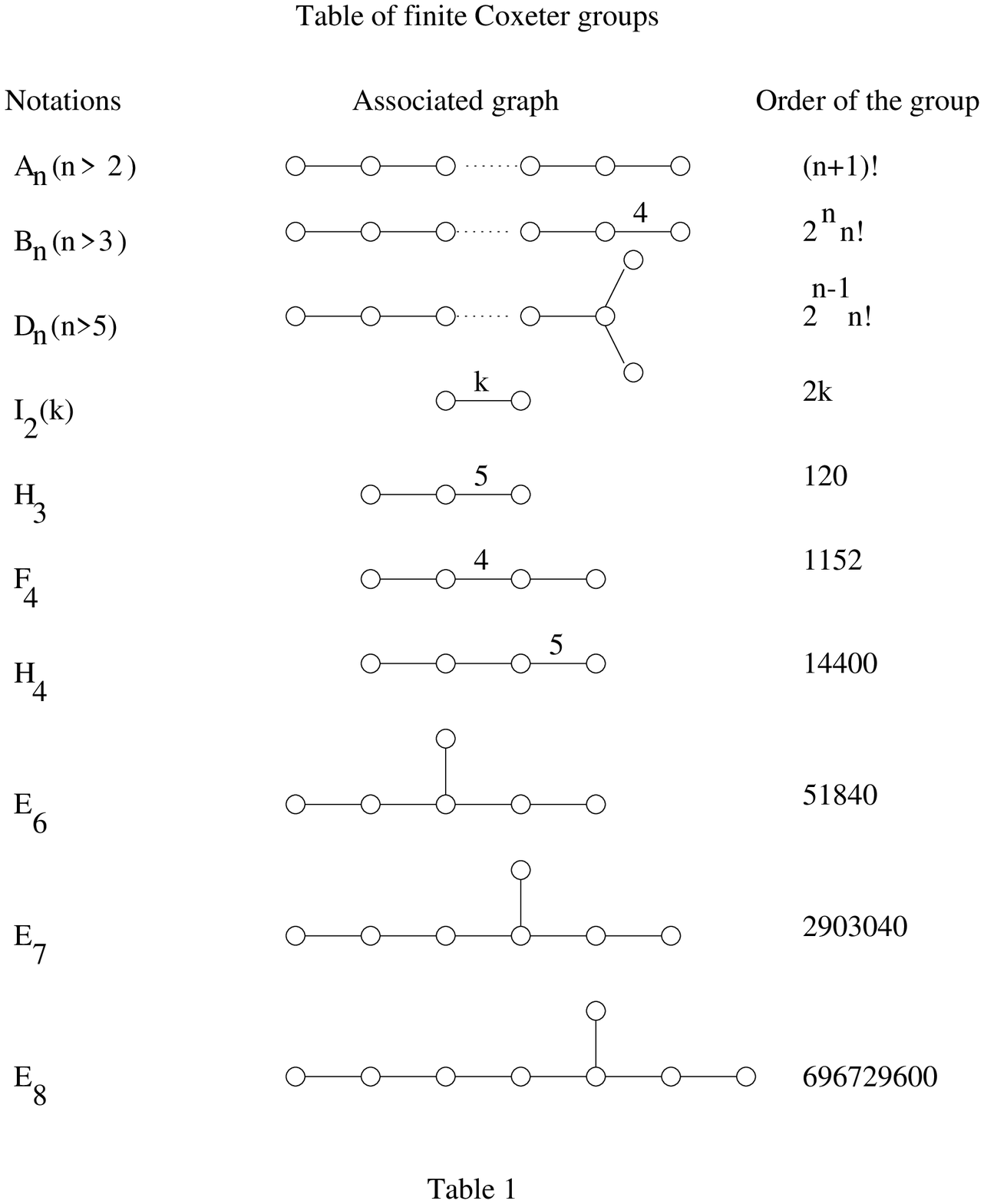,height=15cm,width=12cm}}

\vskip 3cm

It is a long standing question if the classical Artin braid group has 
a discrete faithful linear representation. R. 
Rouquier has informed the author that recently S. Bigelow has proved that
the Artin braid group is linear using Krammer's representation, a 
representation that arises from the BMW algebra. The
well known Burau representation of the Artin braid groups ${\Cal A}_n$ 
in the matrix group
is known to be not faithful for $n\geq 6$ (for $n\geq 10$ see \cite{10}
and for $n\geq 6$ see \cite{9}). 

Recently M. Bestvina (\cite{2}) asked the following two questions: 
(1) Does ${\Cal
A}_S$ satisfy the Tits alternative? (2) Is ${\Cal A}_S$ virtually
poly-free?

If the map $f_S$ is an injective homomorphism then the following results 
will follow from the fact that the classical pure braid group is
poly-free:

\proclaim{Consequence 1.3} Any finite type generalized braid group is 
virtually poly-free.\endproclaim

\proclaim{Consequence 1.4} Finite type generalized braid groups satisfy 
Tits alternative. That is if $H\subset {\Cal A}_S$ is a subgroup which is
not virtually abelian, then $H$ necessarily contains a nonabelian free 
group.\endproclaim

Here we recall that the main result in \cite{8} was proved using 
the fact that the classical pure braid group is strongly poly-free.
Hence we ask:

\proclaim{Question 2} Is any generalized braid group virtually 
strongly poly-free?\endproclaim

For curious reader, here we recall the definition of strongly poly-free 
groups:

\proclaim{Definition} {\rm A discrete group $\Gamma$ is called 
{\it strongly poly-free} if there exists a finite filtration of $\Gamma$ 
by subgroups: 
$1=\Gamma_0\subset \Gamma_1\subset \cdots \subset \Gamma_n=\Gamma$ 
such that the following conditions are satisfied:
\roster
\item $\Gamma_i$ is normal in $\Gamma$ for each $i$

\item $\Gamma_{i+1}/\Gamma_i$ is a finitely generated free group 

\item for each $\gamma\in \Gamma$ and $i$ there is a compact 
surface $F$ and a diffeomorphism $f:F\to F$ such that the induced 
homomorphism $f_{\#}$ on $\pi_1(F)$ is equal to $c_\gamma$ in 
$Out(\pi_1(F))$, where $c_\gamma$ is the action of $\gamma$ on 
$\Gamma_{i+1}/\Gamma_i$ by conjugation and $\pi_1(F)$ is
identified with $\Gamma_{i+1}/\Gamma_i$ via a suitable isomorphism.
\endroster
\endproclaim

\remark{Remark 1.5} R. Rouquier has informed the author about their 
paper \cite{5} on finite complex reflection groups. There they have 
shown, among other results, that except for few cases, for finite complex
reflection groups 
there are presentation of the groups in terms of generators and relations 
so that the associated generalized braid group can be presented after 
removing only the finite order relations on the generators in the
reflection group presentation. For example this is true in the case of
real reflection groups (i.e., Coxeter groups). 
So the same Questions 1 and 2 may be asked for the generalized braid
groups associated to these finite complex reflection groups.

Here recall that a finite 
reflection subgroup $G\subset GL(n, {\Bbb C})$ has a finite set of
generators $\{s\ |\ s\in S\}$ so that each $s$ fixes (pointwise) a
hyperplane $H_s$ in ${\Bbb C}^n$. The group $G$ acts on the space 
${\Bbb C}^n-\cup_{s\in S}H_s$ fixed point freely. Consider the quotient 
space $({\Bbb C}^n-\cup_{s\in S}H_s)/G$. The fundamental group of the 
space $({\Bbb C}^n-\cup_{s\in S}H_s)/G$ is by definition the 
generalized braid group associated to the reflection group $G$. It is 
not yet known, if the generalized braid group corresponding to any
finite reflection group $G$ can be obtained from $G$ after removing some
of the relations in $G$. See \cite{5} for details on this subject. 
\endremark

\remark{Remark 1.6} As far as I know the results in the Consequences
1.1-1.4 are still open.\endremark

\head
2. $f_S$ is a homomorphism for the braid groups of type $A_n$, $I_2(k)$ 
\endhead

The Coxeter groups of type $A_n$ are the symmetric groups. So there is 
nothing to prove in this case.

Before we start with the proof of the Theorem in the case $I_2(k)$ we
recall (with a different notation) 
the presentation of 
the classical Artin braid group ${\Cal A}_n$ on $(n-1)$-strings. 
It is generated by 
the symbols $\sm_1, \sigma_2,\cdots, \sigma_{n-1}$ with the 
relations: $\sigma_i\sigma_j=\sigma_j\sigma_i$ for $|i-j|\geq 2$ and 
$\sigma_i\sigma_{i+1}\sigma_i=\sigma_{i+1}\sigma_i\sigma_{i+1}$ for 
$i\leq n-2$. Geometrically the generator $\sigma_i$ is represented 
by the following braids (Figure 1):

\medskip

\centerline{\psfig{figure=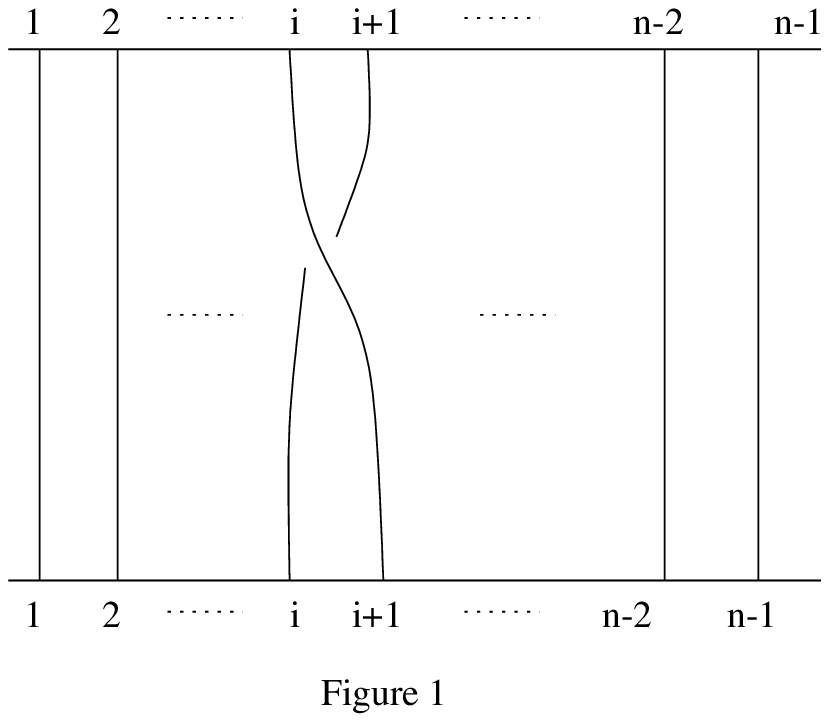,height=6cm,width=7cm}}

\vskip 1cm

The multiplication is given by juxtaposition of braids. And two 
braids are said to be same if one can be deformed to the other 
by an isotopy (fixing the end points $1,2,\cdots, n-1$) without any 
crossing of any two threads of the braids (see Chapter 1 in \cite{3} 
for more details on this.) 
The surjective homomorphism ${\Cal A}_n\to S_n$ is defined by sending 
$\sigma_i$ to the transposition $(i,i+1)$.
 
Now we recall the presentation of the Coxeter groups of type $I_2(k)$: 
${\Cal W}_S=\{s_1, s_2\ |\ s_1^2=s_2^2=(s_1s_2)^k=1\}$. First we consider 
the case when $k$ is even. The case of $k$ odd will be treated later 
on.

Note that $|{\Cal W}_S|=2k$. We now consider $k=2$. In terms of 
coset representative we write 
${\Cal W}_S=\{1, s_1, s_2, s_1s_2\}$. Consider the bijection 
$\{1, 2, 3, 4\}\to \{1, s_1, s_2, s_1s_2\}$ given by $1\mapsto 1, 2\mapsto 
s_1, 3\mapsto s_2, 4\mapsto s_1s_2$. It is easy to check that there 
is an embedding ${\Cal W}_S\to S_4$ sending $s_1\mapsto (12)(34)$ and 
$s_2\mapsto (13)(24)=(23)(12)(34)(23)$. Define the map  
$f_S:{\Cal A}_S\to {\Cal A}_4$ by $$s_1\mapsto \sigma_1\sigma_3\ 
\text{and}\  
s_2\mapsto \sigma_2\sigma_1\sigma_3\sigma_2$$

To show that 
$f_S$ is a homomorphism we only have to check that the two 
braids $\sigma_1\sigma_3\sigma_2\sigma_1\sigma_3\sigma_2$ and 
$\sigma_2\sigma_1\sigma_3\sigma_2\sigma_1\sigma_3$ are same. 
We draw the pictures (Figure 2) of these two braids and see that they are 
same upto isotopy.

Now we consider the case of $k=4$. In this case the group ${\Cal W}_S$ 
has $8$ elements. We give an embedding of this  group to the symmetric 
group on $8$ letters. We write in terms of coset representative ${\Cal
W}_S=\{1,s_1,s_2,s_1s_2,s_2s_1,
s_1s_2s_1,s_2s_1s_2,s_1s_2s_1s_2\}$ and the bijection between 
${\Cal W}_S$ 
and $\{1,2,3,4,5,6,7,8\}$ preserving the order in the way the elements 
are written. The embedding ${\Cal W}_S\to S_8$ is given by 
$s_1\mapsto (12)(34)(56)(78)$ and $$s_2\mapsto (13)(25)(47)(68)=
(23)(12)(45)(34)(23)(67)(56)(45)(78)(67)$$

Define the map 
$f_S:{\Cal A}_S\to {\Cal A}_8$ by $s_1\mapsto
\sigma_1\sigma_3\sigma_5\sigma_7$ and $$s_2\mapsto
\sigma_2\sigma_1\sigma_4
\sigma_3\sigma_2\sigma_6\sigma_5\sigma_4\sigma_7\sigma_6$$
 
To check that $f_S$ is a homomorphism we need to show that
$f_S(s_1s_2s_1s_2)=f_S(s_2s_1s_2s_1)$. Again we actually show by 
drawing pictures (Figures 3,4,5) that the two associated braids are 
same upto isotopy.

The case of $k=6$ is similar. We just write down the embedding of 
${\Cal W}_S$ in $S_{12}$: $s_1\mapsto (12)(34)(56)(78)(910)(1112)$ and 

$$s_2\mapsto (13)(25)(47)(69)(811)(1012)=$$
$$(23)(12)(45)(34)(23)(67)(89)(78)(67)(1011)(910)(89)(1112)(1011)$$ 

We define $f_S:{\Cal A}_S\to {\Cal A}_{12}$ by $f_S(s_1)=\sm_1\sm_3\sm_5
\sm_7\sm_9\sm_{11}$ and 

$$f_S(s_2)=\sm_2\sm_1\sm_4\sm_3\sm_2\sm_6
\sm_8\sm_7\sm_6\sm_{10}\sm_9\sm_8\sm_{11}\sm_{10}$$ 

Again we need to 
check that $f_S(s_1s_2s_1s_2s_1s_2)=f_S(s_2s_1s_2s_1s_2s_1)$. We show 
it in Figures 6 and 7. 

In the case of $k=2$ it is easy to see (by drawing the picture of
the braid $f_S(s_1^ms_2^n)$) that the map $f_S$ is in fact injective. See 
Figure 8.

Now we write down the embedding $e$ for any even $k$: 
$e:{\Cal W}_S\to S_{2k}$ is defined by $$s_1\mapsto (12)(34)(56)(78)\cdots 
(2k-1, 2k)$$ and $$s_2\mapsto (13)(25)(47)(69)\cdots (2k-4, 2k-1)
(2k-2, 2k)$$

Define $f_S$ by: $f_S(s_1)=\sm_1\sm_3\ldots\sm_{2k-1}$ and 
to define $f_S(s_2)$ we do not need to represent it in terms of 
$\sm_i$, but at first we draw the lines from $1$ to $3$, $2$ to $5$, 
$4$ to $7$,...,$2k-2$ to $2k$ then we draw the lines from $3$ to $1$, 
$5$ to $2$,...,$2k$ to $2k-2$ and pass through below if any previous
line comes on the way. 

It can easily be seen after drawing the pictures of 
the braids $f_S(s_1s_2\ldots s_1s_2)$ and $f_S(s_2s_1\ldots s_2s_1)$ 
that they are same up to isotopy. In fact they are of the type as the 
second picture in Figure 7. 

Finally we consider the case of when $k$ is odd. In this case 
we choose the embedding $e$ in the following way: 
$e:{\Cal W}_S\to S_k$ is defined by $$s_1\mapsto (23)(45)(67)\cdots 
(k-1, k)$$ and $$s_2\mapsto (12)(34)(56)\cdots (k-2, k-1)$$ 

We define the map $f_S$ by: $f_S(s_1)=\sm_2\sm_4\cdots\sm_{k-1}$ 
and $f_S(s_2)=\sm_1\sm_2\cdots\sm_{k-2}$. 

We draw the pictures of the braids $f_S(s_1s_2\ldots s_1s_2s_1)$ and 
$f_S(s_2s_1\ldots s_2s_1s_2)$ for some special cases in Figure 9 and 10 
and see that they are same upto isotopy. The general case is easily 
seen to be true. 

Note that $k=3$ is the symmetric groups case. So we consider the 
special cases $k=5,7$. In these cases the map $f_S$ is of the following 
type: $k=5$: $$f_S(s_1)=\sm_2\sm_4\ \text{and}\ f_S(s_2)=\sm_1\sm_3$$ and 
$k=7$: $$f_S(s_1)=\sm_2\sm_4\sm_6\ \text{and}\ f_S(s_2)=\sm_1\sm_3\sm_5$$ 

Figure 9 shows the case $k=5$ and Figure 10 corresponds to $k=7$.

\newpage

\centerline{\psfig{figure=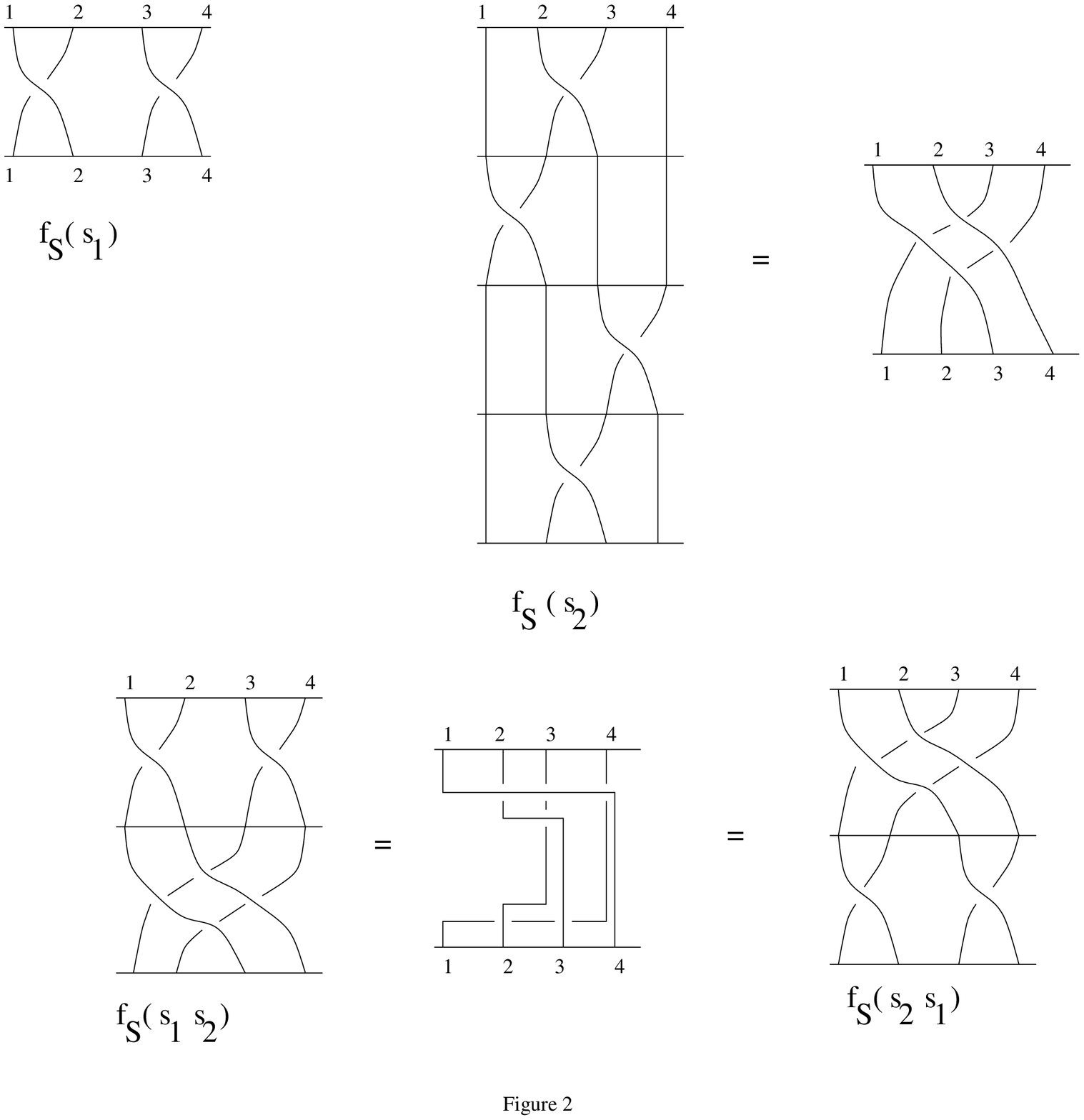,height=17cm,width=12.5cm}}

\newpage

\centerline{\psfig{figure=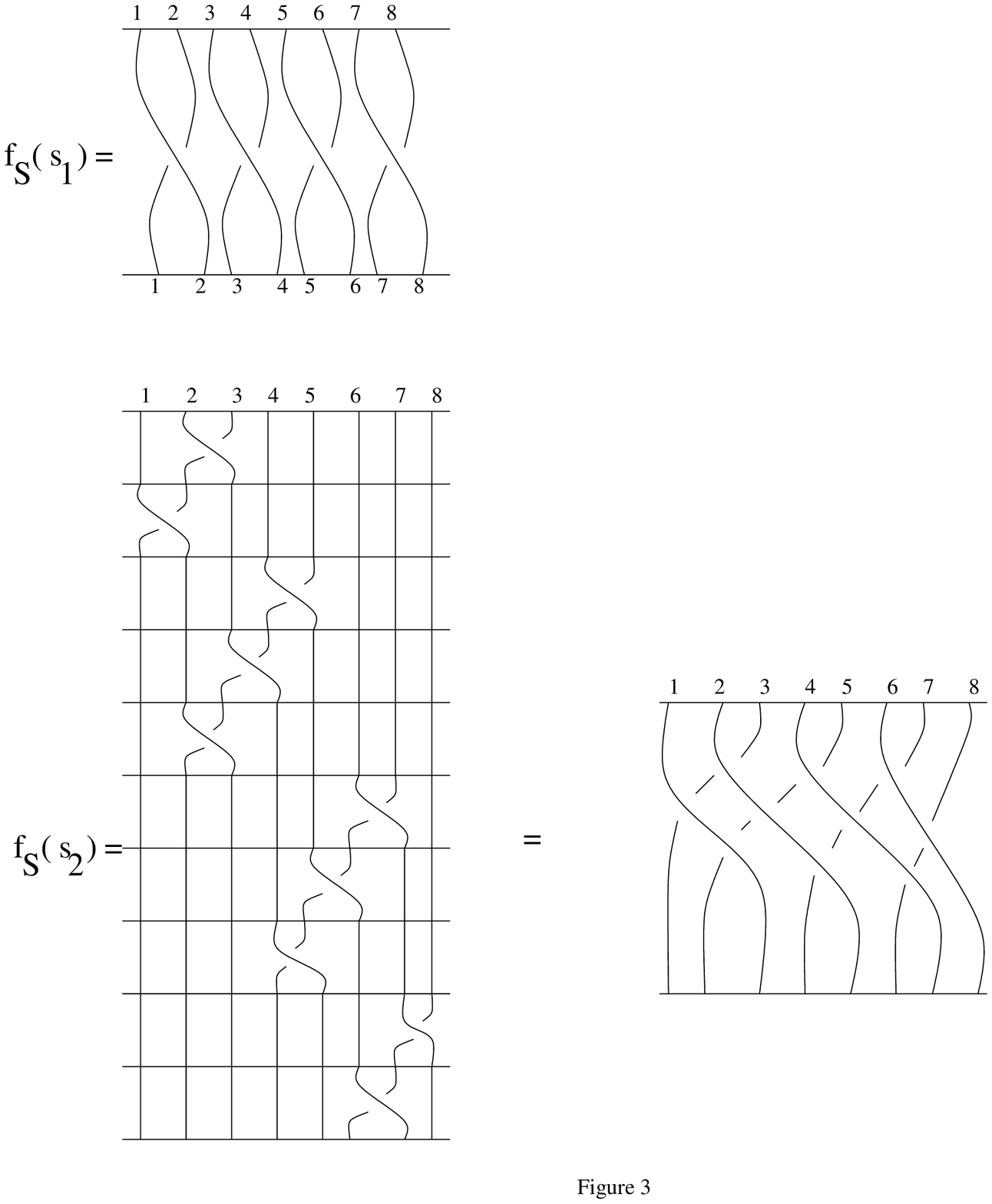,height=18cm,width=12.5cm}}

\newpage

\centerline{\psfig{figure=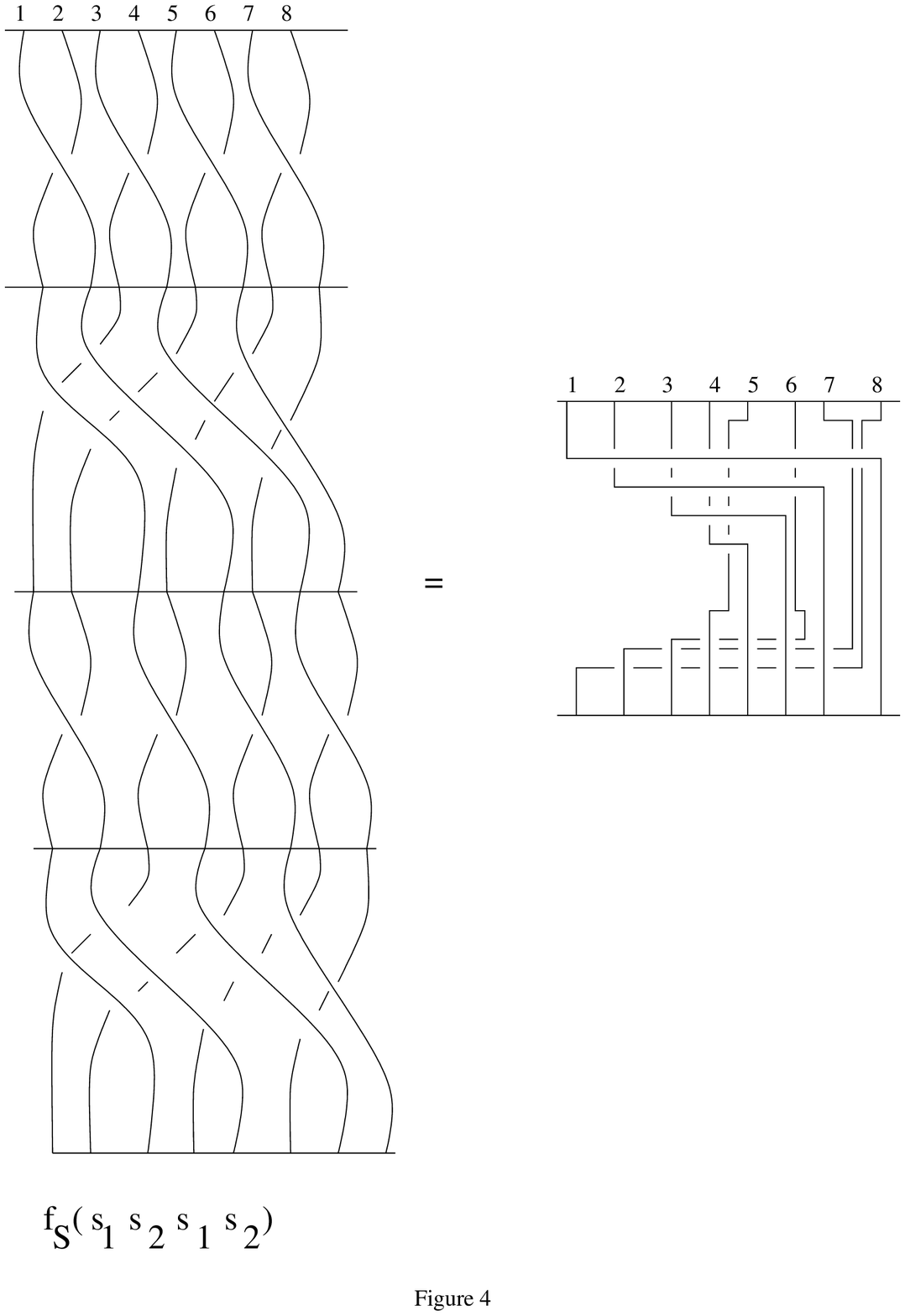,height=18cm,width=12.5cm}}

\newpage

\centerline{\psfig{figure=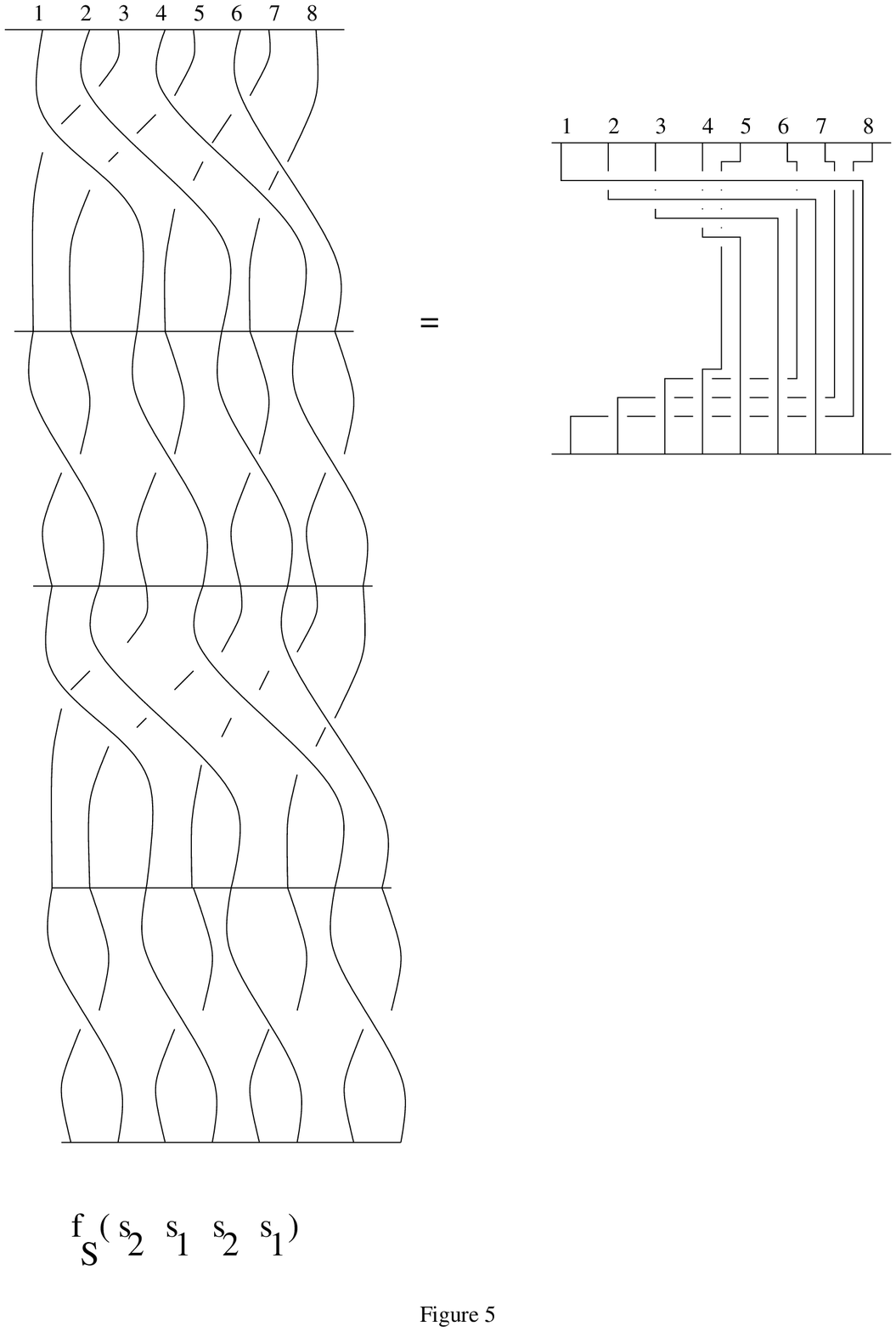,height=18cm,width=12.5cm}}

\newpage

\centerline{\psfig{figure=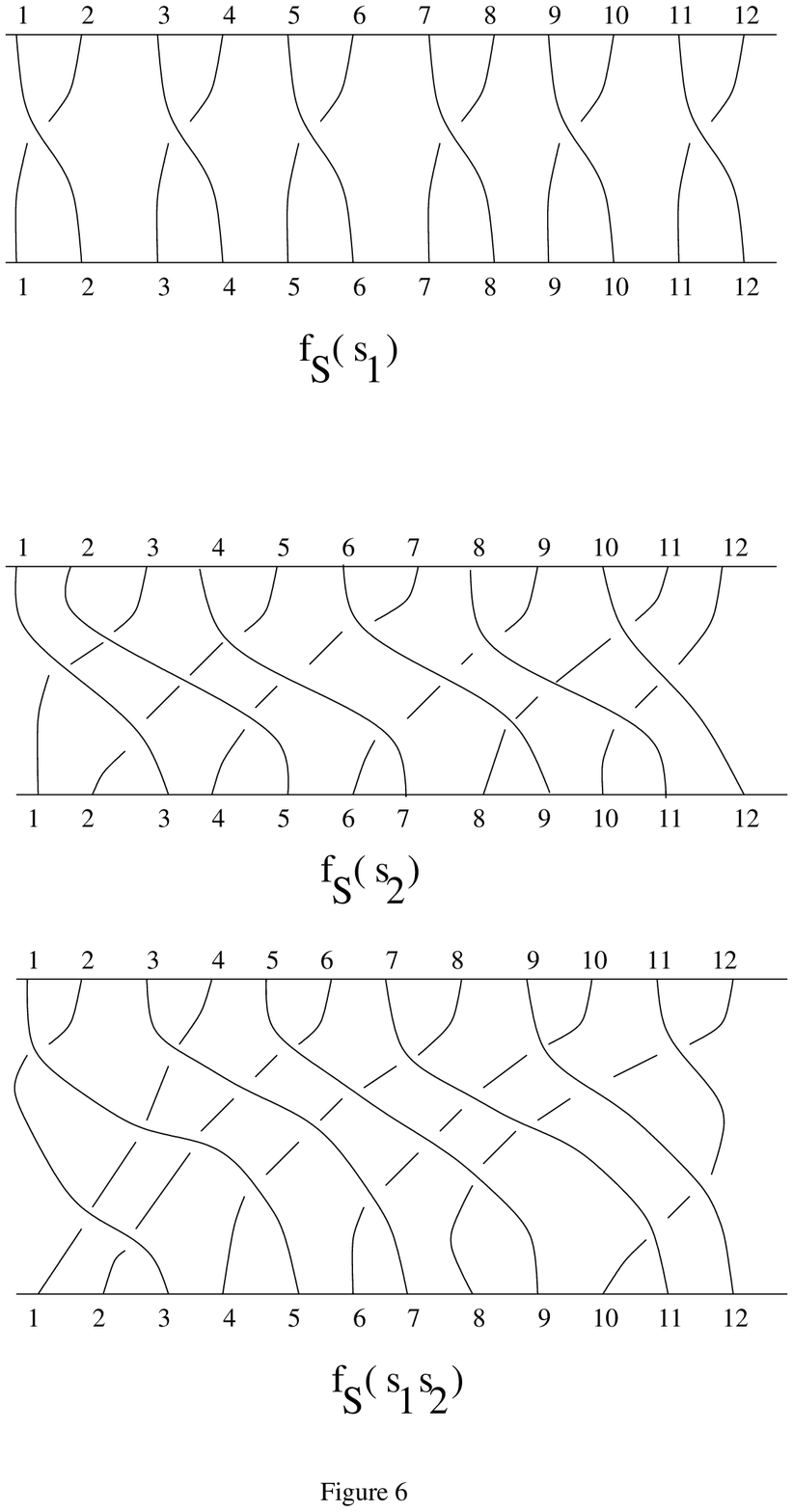,height=18cm,width=12.5cm}}

\newpage

\centerline{\psfig{figure=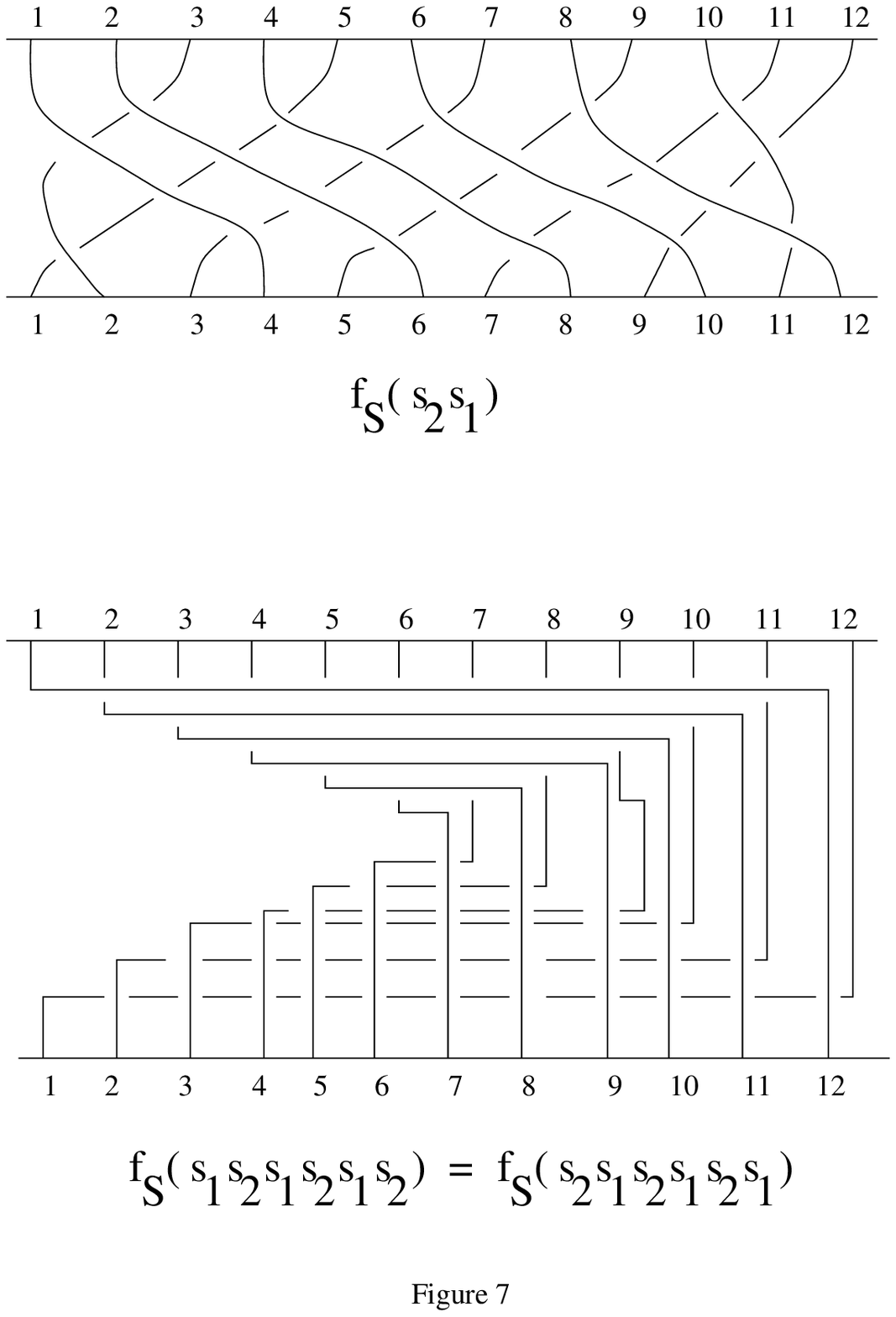,height=18cm,width=12.5cm}}

\newpage

\centerline{\psfig{figure=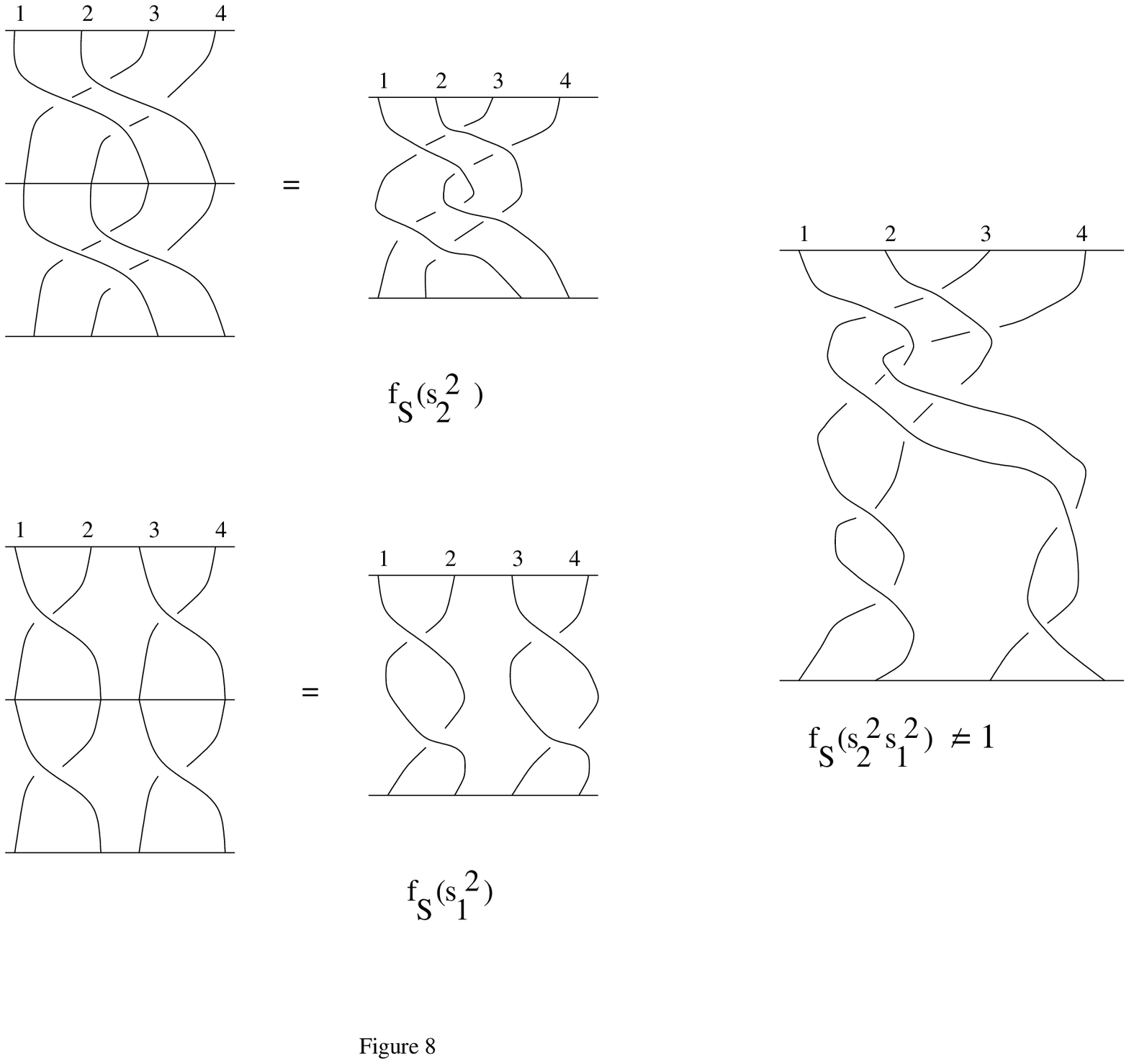,height=16cm,width=12.5cm}}

\newpage

\centerline{\psfig{figure=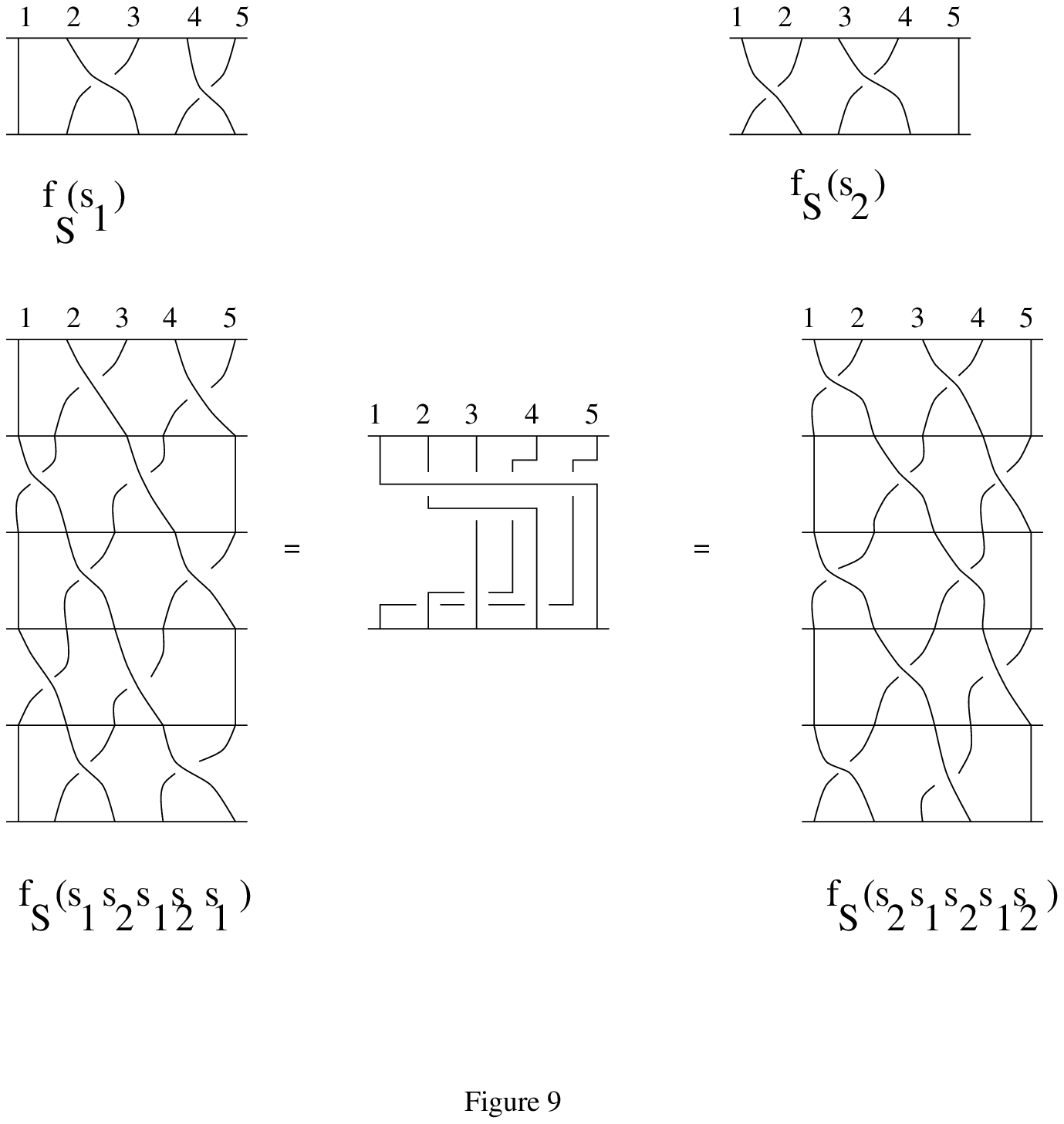,height=16cm,width=12.5cm}}

\newpage

\centerline{\psfig{figure=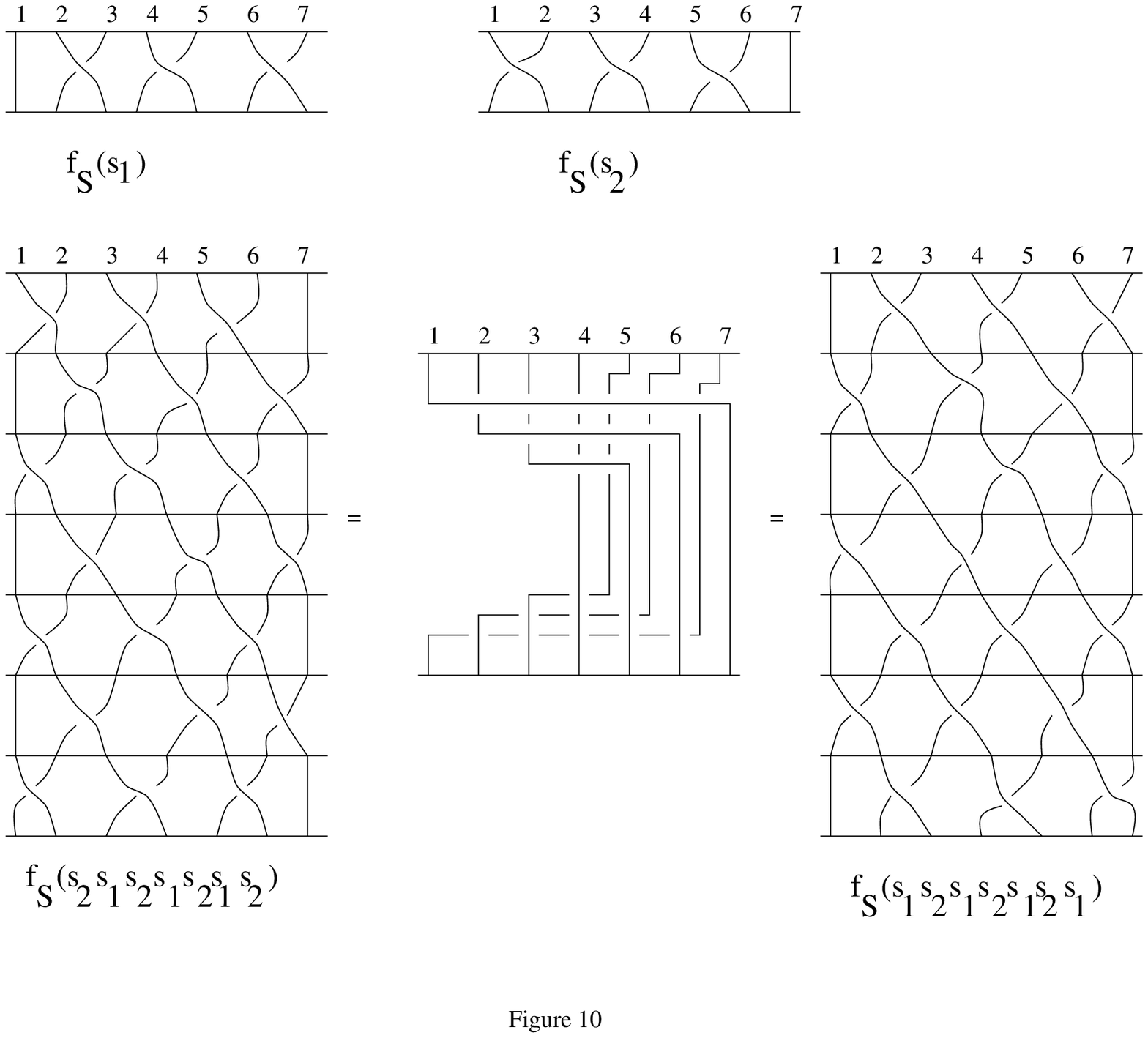,height=15cm,width=12.5cm}}

\newpage

Hence we have proved the following theorem:

\proclaim{Theorem 2.1} Let ${\Cal W}_S=\{s_1,s_2\ |\ s_1^2=
s_2^2=(s_1s_2)^k=1\}$ be the Coxeter group of type $I_2(k)$. 
Then there is a homomorphism $f_S:{\Cal A}_S\to {\Cal A}_m$ and an 
embedding $e:{\Cal W}_S\to {\Cal S}_m$, for 
some $m$ (depending on $k$), making the diagram $D$ 
commutative. Also the map $f_S$ is injective in the case $k=2$.
\endproclaim  

\head
3. $f_S$ is a homomorphism for the braid groups of type $B_n$
\endhead

Suppose ${\Cal W}_S$ be a Coxeter group of type $B_n$. It has a 
representation of the following form: 
$${\Cal W}_S=\{s_1, s_2, \cdots, s_n\
|\ s_1^2=s_2^2=\cdots= s_n^2=(s_1s_2)^4=(s_is_{i+1})^3\  
(\text{for}\ i\geq 2)=$$ $$(s_is_j)^2\ (\text{for}\ |i-j|\geq 2)=1\}.$$ 

There is a faithful representation $e:{\Cal W}_S\to S_{2n}$ defined 
on the generators in the following way: $s_{n-j}\mapsto (j+1, j+2)(2n-j-1,
2n-j)$ for $2\leq n-j\leq n$ and $s_1\mapsto (n, n+1)$. 

We define the map $f_S:{\Cal A}_S\to {\Cal A}_{2n}$ by $f_S(s_{n-j})=
\sm_{j+1}\sm_{2n-j-1}$ and $f_S(s_1)=\sm_n$. 

We check that $f_S$ is a homomorphism in Figure 11. We show the 
following equivalence of braids: $f_S(s_1s_2s_1s_2)=f_S(s_2s_1s_2s_1)$, 
$f_S(s_is_{i+1}s_i)=f_S(s_{i+1}s_is_{i+1})$ for $i\geq 2$ and 
$f_S(s_is_j)=f_S(s_js_i)$ for $|i-j|\geq 2$. The first two equivalence 
follows from Figure 11 and the last equivalence follows from the 
braid relation $\sm_i\sm_j=\sm_j\sm_i$ for $|i-j|\geq 2$.

Thus we have proved the following theorem:

\proclaim{Theorem 3.1} Let ${\Cal W}_S$ be a Coxeter group 
of type $B_n$. Then there is a homomorphism 
$f_S:{\Cal A}_S\to {\Cal A}_{2n}$ and an
embedding $e:{\Cal W}_S\to {\Cal S}_{2n}$, making the  
diagram $D$ commutative. 
\endproclaim

\newpage

\centerline{\psfig{figure=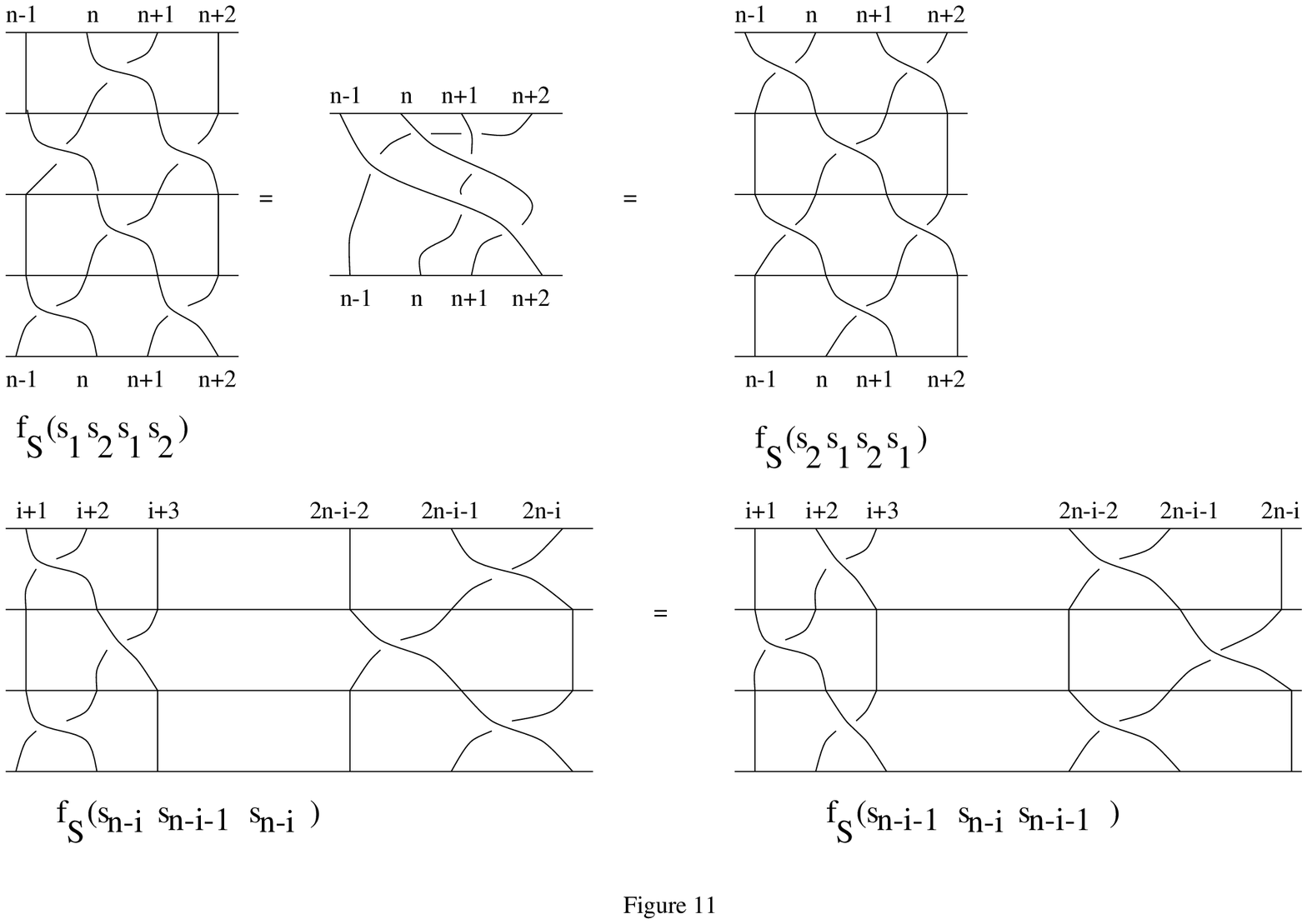,height=15cm,width=14cm}}

\newpage

\head
4. An example of an extension $f_S$ which is not a homomorphism
\endhead

Consider the Coxeter group of type $D_n$. It has a representation 
of the following type: $${\Cal W}_S=\{s_1,\ldots, s_n\ | \ (s_1s_2)^3=
(s_2s_3)^3=(s_2s_4)^3=(s_is_{i+1})^3 (\text{for}\ i\geq 4)=$$ $$(s_1s_i)^2
(\text{for}\ i\geq 3)=(s_3s_i)^2 (\text{for}\ i\geq 4)=(s_is_j)^2
(\text{for}\ |i-j|\geq 2\ \text{and}\ i,j\geq 4)=1\}.$$ We choose the case
$n=4$. In this case there is an embedding of $D_4$ in the symmetric 
group $S_8$ given by $s_1\mapsto (34)(56), s_2\mapsto (23)(67), 
s_3\mapsto (35)(46), s_4\mapsto (12)(78)$. Note that
$(35)(46)=(45)(34)(56)(45)$. We define $f_S$ by $f_S(s_1)=\sm_3\sm_5, 
f_S(s_2)=\sm_2\sm_6, f_S(s_3)=\sm_4\sm_3\sm_5\sm_4, f_S(s_4)=\sm_1\sm_7$. 
In Figure 12 we check that the two braids $f_S(s_2s_3s_2)$ and
$f_S(s_3s_2s_3)$ are not same. It is easy to see in Figure 12 that 
the strings $3$ to $3$ and $6$ to $6$ in $f_S(s_2s_3s_2)$ are not 
tangled, whereas the strings $3$ to $3$ and $6$ to $6$ in the braid 
$f_S(s_3s_2s_3)$ are tangled which cannot be untangled by an isotopy 
of braids. 

Thus we have showed that the natural map $f_S$ is not a homomorphism, 
but the diagram $D$ commutes. At the time of writing this paper  
it was not clear how to choose this extension to be a homomorphism 
so that the diagram $D$ commute. 

\remark{Remark 4.1} For the generalized braid groups corresponding to the 
Coxeter groups of type $H_3$, $H_4$, $F_4$, $E_6$, $E_7$ and $E_8$ the
author found some embedding of the Coxeter group in some symmetric group
but either the symmetric group was too large to draw the pictures of the
braids for checking if an extension $f_S$ is a homomorphism or a similar 
situation as in the case of $D_n$ (as above) occurred. Also we
note that $H_3$ is a parabolic subgroup of $H_4$, thus it is
enough to check if $f_S$ is an injective homomorphism in the
case of $H_4$ to prove the same for $H_3$.\endremark

\remark{Remark 4.2} We used the $GAP$ programme (\cite{12}) to find 
the embedding $e:{\Cal W}_S\to S_m$ in some concrete cases of the 
Coxeter group ${\Cal W}_S$ and this was used to find the  
embedding in general for the Coxeter groups $I_2(k)$ (for $k$ odd) and
$B_n$.\endremark

\newpage

\centerline{\psfig{figure=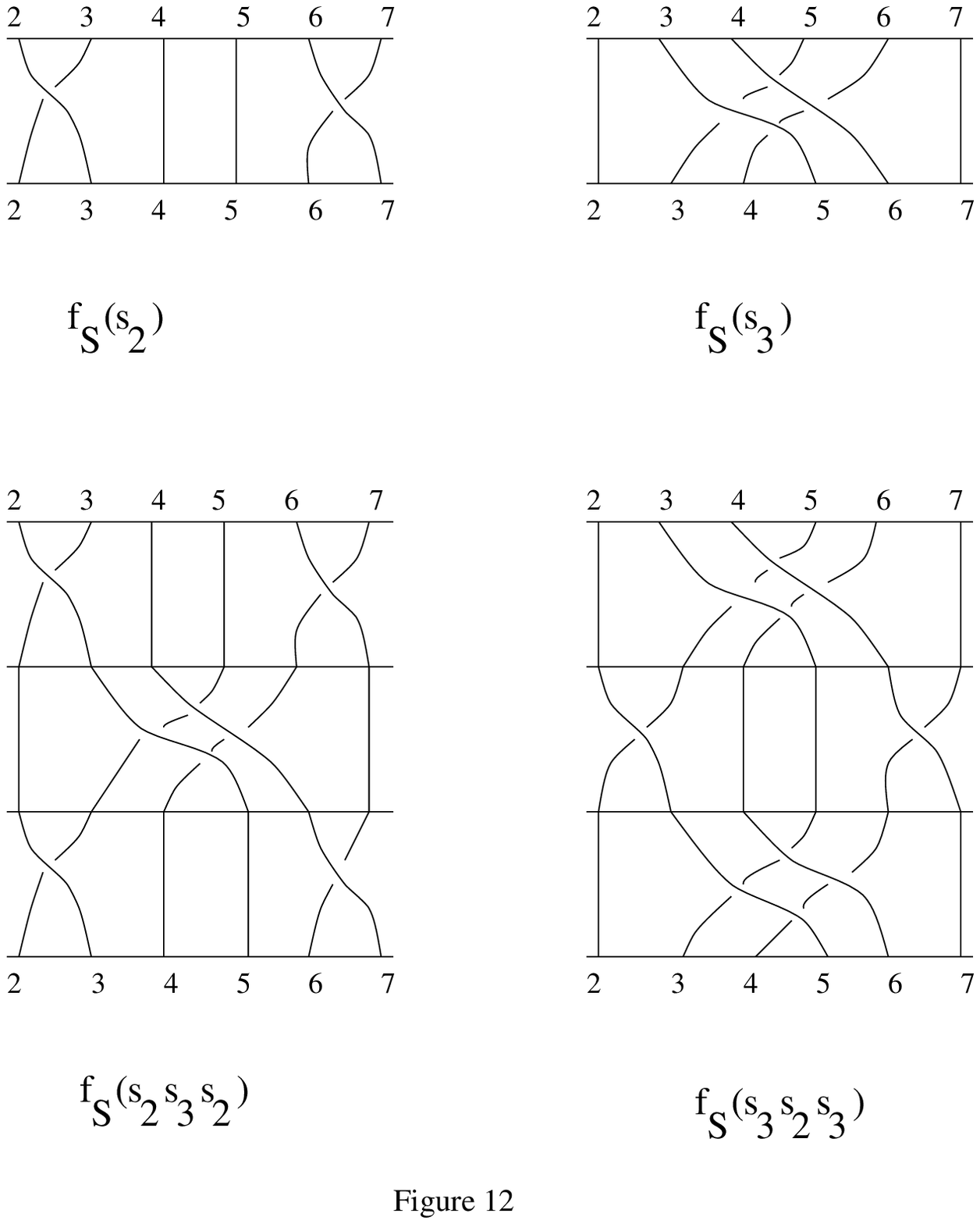,height=16cm,width=14cm}}

\newpage

\remark{Acknowledgment} The author would like to thank the Institut Des
Hautes \'{E}tudes Scientifiques in France for the kind invitation, where
this work was done.\endremark

\enddocument